\documentclass[11pt]{amsart}

\usepackage{amsmath, amsthm, amssymb, amscd, amsfonts}
\usepackage[all]{xy}
\usepackage{hhline}
\usepackage[dvips, dvipsnames, usenames]{color}

\usepackage[T1]{fontenc} 

\usepackage[active]{srcltx}

\newcommand{\supp}{\operatorname{supp}}

\newcommand{\GL}{\operatorname{GL}}
\newcommand{\SL}{\operatorname{SL}}
\newcommand{\PGL}{\operatorname{PGL}}
\newcommand{\PSL}{\operatorname{PSL}}

\newcommand\toba{{\mathfrak B }}

\newcommand{\trid}{\triangleright}

\newcommand{\ku}{\mathbb C}

\newcommand{\Z}{{\mathbb Z}}
\newcommand{\N}{{\mathbb N}}

\newcommand{\G}{{\mathbb G}}

\newcommand{\F}{{\mathbb F}}
\newcommand{\C}{{\mathcal C}}

\newcommand{\oc}{{\mathcal O}}
\newcommand{\ocs}{{\mathcal O}_{\sigma}}

\newcommand{\la}{\langle}
\newcommand{\ra}{\rangle}

\theoremstyle{plain}

\newtheorem{lema}{Lemma}[section]
\newtheorem{theorem}[lema]{Theorem}
\newtheorem{cor}[lema]{Corollary}

\theoremstyle{definition}

\theoremstyle{remark}
\newtheorem{obs}[lema]{Remark}
\newtheorem{rmk}[lema]{Remarks}

\newcommand\am{\mathbb A_m}

\newcommand\A{\mathbb A}

\newcommand\s{\mathbb S}

\def\pf{\begin{proof}}
\def\epf{\end{proof}}

\theoremstyle{remark}

\begin{document}

\renewcommand{\baselinestretch}{1.2}

\thispagestyle{empty}

\title[Conjugacy classes of $p$-cycles of type D]{Conjugacy classes of $p$-cycles of type D\\ in alternating groups}

\author[Fernando Fantino]{Fernando Fantino}

\thanks{This work was partially supported by CONICET, SeCyT - Universidad Nacional de C\'ordoba, MinCyT, Embassy of France in Argentina and Mairie de Paris}

\address{D\'epartement de Math\'ematiques, Universit\'e Paris Diderot (Paris 7),
175, rue du Chevaleret, 75013, Paris, France}


\address{Facultad de
Matem\'atica, Astronom\'{\i}a y F\'{\i}sica, Universidad Nacional
de C\'ordoba, CIEM -- CONICET. Medina Allende s/n (5000) Ciudad
Universitaria, C\'ordoba, Argentina}
\email{fantino@famaf.unc.edu.ar}


\begin{abstract}
We classify the conjugacy classes of $p$-cycles of type D in alternating groups. This finishes the open cases in \cite{AFGV-alt}.
Also we determine all the subracks of those conjugacy classes which are not of type D.
\end{abstract}
\maketitle


\section{Introduction}\label {se:introduction}

In the context of the Lifting method \cite{AS-camb}, the problem of the classification of finite-dimensional complex pointed Hopf algebras over non-abelian groups can be
approached by the study of Nichols algebras associated to pairs $(X,q)$, where $X$ is a rack and $q$ a 2-cocyle, see \cite{AG}.
Since the computation of all cocycles of a rack is hard, it is useful to have tools that ensure that
the corresponding Nichols algebra $\toba(X,q)$ has infinite dimension for any $2$-cocycle $q$; we say that $X$ \emph{collapses} if this happens.
It was shown in \cite[Thm. 3.6]{AFGV-alt} that any finite rack of \emph{type D} collapses.

The racks of type D have a nice behavior with respect to monomorphisms and
epimorphisms of racks. Indeed, if $Y \subseteq X$ is a subrack of type D, then $X$ is of type D, and
if $p:Z\to X$ is an epimorphism of racks with $Z$ finite and $X$ of type D, then $Z$ is of type D.
Besides, it is well-known that any finite rack can be decomposed as a union of indecomposable subracks and that
every indecomposable rack $X$ admits a projection $X\to Y$ with $Y$
a \emph{simple} rack, i.~e. a rack without proper quotients. We recall also that the classification of finite simple racks is known, see
\cite{AG} and \cite{Jo}.

These facts suggest that the notion of racks of type D is useful
for an approach for the classification problem of
finite-dimensional pointed Hopf algebras over non-abelian groups and it establish a first step for that problem: to classify
finite simple racks of type D, see \cite[\S 2.6]{AFGaV-cla}.

One of the most important family of finite simple racks are the conjugacy classes of finite non-abelian (almost) simple groups. In \cite{AFGV-alt}, all the conjugacy classes of type D in symmetric and alternating groups were determined, except the conjugacy classes of $p$-cycles in $\A_p$ and $\A_{p+1}$.

This article is a contribution to the problem of classification of finite simple racks of type D.
Explicitly, we determine when the conjugacy classes of elements of order $p$, with $p\geq 5$ prime, in the alternating groups $\A_p$ and $\A_{p+1}$ are
of type D or not.

We summarize our results in the following statement.

\begin{theorem}\label{theorem:main}
Let $p$ be a prime number, $p\geq 5$, and $m\in\{p,p+1\}$. Let $\oc$ be a conjugacy class of $p$-cycles in $\A_m$.
\begin{enumerate}
\renewcommand{\theenumi}{\Roman{enumi}}\renewcommand{\labelenumi}{(\theenumi)}
  \item If $m=p$, then $\oc$ is of type D if and only if $p\geq 13$ and $p=\frac{r^k-1}{r-1}$, with $r$ a prime power and $k$ is a natural number.
  \item If $m=p+1$, then $\oc$ is of type D if and only if $p\geq 7$ and $p=\frac{r^k-1}{r-1}$, with $r$ a prime power and $k$ is a natural number.
\end{enumerate}
\end{theorem}

This result finishes the open cases in \cite[Thm. 4.1]{AFGV-alt}, i.~e. it concludes the classification of conjugacy classes of type D in alternating groups.

The family of primes $p$ of the form $\frac{r^k-1}{r-1}$, with $r$ a power of a prime number, contains the Mersenne primes and the Fermat primes.
The primes $p$ of this form with $p<1000$ are: 3, 5, 7, 13, 17, 31, 73, 127, 257, 307 and 757, see Remark \ref{obs:ciclot:primos:1000}.

\section{Preliminaries}\label {se:preliminaries}

Throughout the paper $M_{11}$, $M_{12}$, $M_{23}$ and $M_{24}$ denote the corresponding Mathieu
simple groups and $L_k(r)$ means the projective special linear group, $r$ a prime power.
For $m\in \N$, $\G_m$ denotes the $m$-th roots of 1 in $\ku$.

\subsection{}

A \emph{rack} is a pair $(X,\trid)$, where $X$ is a non-empty set
and $\trid:X\times X \to X$ is a function, satisfying the
following conditions: for every $x\in X$, the function $x\trid -
:X \to X$ is bijective and $x\trid (y\trid z)=(x\trid y) \trid
(x\trid z)$, for all $x$, $y$, $z\in X$.
Any subset of a group $G$ stable by conjugation is a rack with the
conjugation as function $\trid$. In particular, a conjugacy class
of $G$ is a rack.

A rack $(X,\trid)$ is
said to be \emph{of type D}\footnote{The letter D stands for decomposable.} if it contains a decomposable subrack
$Y = R\coprod S$ such that $r\trid(s\trid(r\trid s)) \neq s$, for
some $r\in R$, $s\in S$.
It is easy to see that a conjugacy class $\oc$
of a group $G$ is a rack of type D if and only if there exist $\sigma$, $\tau\in \oc$
such that
\begin{align*}
&\text{{\bf (Ax.~1)}} & &(\sigma\tau)^2\neq (\tau\sigma)^2, &&&&&&&\\
&\text{{\bf (Ax.~2)}} & &\text{ $\sigma$ and $\tau$ are not conjugated in $\la\sigma,\tau\ra$},&&&&&&&
\end{align*}
where $\la\sigma,\tau\ra$ means the subgroup generated by $\sigma$ and $\tau$.

\medbreak

The importance of studying racks of type D lies on the following result.

\begin{theorem} \cite[Thm. 3.6]{AFGV-alt}
If $X$ is a finite rack of type D, then $\toba(X,q)$ has infinite dimension for all $2$-cocycle $q$. \qed
\end{theorem}

This result is based on \cite[Thm. 8.6]{HS}, a consequence of \cite{AHS}.

\subsection{}
Let $\s_m$ and $\A_m$ be the symmetric group and alternating group
in $m$ letters, respectively. Let $\sigma\in \s_m$. It is well-known that
the conjugacy class $\ocs^{\s_m}$ of $\sigma$ in $\s_m$
coincides with the set of permutations in $\s_m$ with the same \emph{type}, i.~e. the cycle structure, as $\sigma$.
On the other hand, if $\sigma\in\A_m$ and $\ocs^{\A_m}$ denotes the conjugacy class of $\sigma$ in $\A_m$, then either
$\ocs^{\s_m}=\ocs^{\A_m}$ or else $\ocs^{\s_m}=\ocs^{\A_m} \cup
\oc_{(1\, 2)\trid \sigma}^{\A_m}$ a disjoint union of two
conjugacy classes in $\am$. This last case occurs if and only if
$\sigma$ is a product of disjoint cycles whose lengths
are odd and distinct.

Let $m\geq 5$. In \cite[Thm. 4.1]{AFGV-alt}, it was proven that

\begin{itemize}
\item if $\sigma\in\A_m$ and the type of $\sigma$ is different from
\begin{enumerate}
\renewcommand{\theenumi}{\roman{enumi}}\renewcommand{\labelenumi}{(\theenumi)}
\item \label{classesAm} $(3^2)$, $(2^2,3)$, $(1^{m-3},3)$, $(2^4)$, $(1^2,2^2)$, $(1,2^2)$,
\item \label{classesAm:primo} $(p)$, $(1,p)$, $p$ prime,
\end{enumerate}
then  $\ocs^{\A_m}$ is of type D;
\item if $\sigma\in\s_m$ and the type of $\sigma$ is different from (i), (ii) and
\begin{enumerate}
\item[(iii)] $(2,3)$,  $(2^3)$,  $(1^{m-2},2)$,
\end{enumerate}
then  $\ocs^{\s_m}$ is of type D.
\end{itemize}
The classes in (i) and (iii) above are not of type D, see \cite[Rmk. 4.2]{AFGV-alt}.

\medbreak

In the present paper we are concerned about the remaining cases: the conjugacy class of $p$-cycles, $p$ prime, in $\A_m$ and in $\s_m$ with $m\in \{p,p+1\}$. For some values of $p$ the problem was already considered in \cite{AFGV-alt}:
\begin{itemize}
        \item if the type of $\sigma$ is $(p)$, then $\ocs^{\A_p}$ is of type D for $p= 13$, $17$, $31$, and $\ocs^{\A_p}$ is not of type D for $p=5$, $7$, $11$;
        \item if the type of $\sigma$ is $(1,p)$, then $\ocs^{\A_{p+1}}$ is of type D for $p=2^q-1$ a Mersenne prime, and $\ocs^{\A_{p+1}}$ is not of type D for $p= 5$, $11$.
    \end{itemize}

\begin{rmk}\label{rem:OpOp+1}
(a) The two conjugacy classes of $p$-cycles in $\A_p$ (resp. in $\A_{p+1}$) are isomorphic as racks.

(b) If $\sigma$ is a $p$-cycle and $\ocs^{\A_p}$ is of type D, then $\ocs^{\A_{p+1}}$ is of type D.
\end{rmk}

\subsection{Subgroups of $\A_m$ generated by two $p$-cyles, $p$ prime}

Let $m$, $p\in \N$, $p$ odd prime. For $\sigma\in \A_m$ we
define $\supp(\sigma):=\{i\in\{1,\dots,m\} \,:\,\sigma(i)\neq
i\}$, i.~e. $\supp(\sigma)$ is the set of points in
$\{1,\dots,m\}$ moved by $\sigma$.

The main tool to prove Theorem \ref{theorem:main} is the following result.

\begin{theorem}\label{teor:FW}\cite{FW}
Let $\sigma$, $\tau$ two $p$-cycles in $\A_m$, with $m=|\supp(\sigma)\cup
\supp(\tau)|$. Then one of the following must occur:
\begin{enumerate}
\renewcommand{\theenumi}{\roman{enumi}}\renewcommand{\labelenumi}{(\theenumi)}
    \item \label{FW:i} $m=p$ and $\langle \sigma, \tau \rangle\simeq \Z/p\Z$;
    \item \label{FW:ii} $m=2p$ and $\langle \sigma, \tau \rangle \simeq \Z/p\Z\times
    \Z/p\Z$;
    \item \label{FW:iii} $m=p=\frac{r^k-1}{r-1}$ and $\langle \sigma, \tau
    \rangle\simeq L_k(r)$;
    \item \label{FW:iv} $m=p+2$, $p$ is a Mersenne prime and $\langle \sigma, \tau
    \rangle\simeq L_2(p+1)$;
    \item \label{FW:v} $m=p+1$ and $\langle \sigma, \tau \rangle\simeq
    L_2(p)$;
    \item \label{FW:vi} $m=p=11$ and $\langle \sigma, \tau \rangle\simeq L_2(11)$,
    $M_{11}$;
    \item \label{FW:vii} $m=p=23$ and $\langle \sigma, \tau \rangle\simeq
    M_{23}$;
    \item \label{FW:viii} $m=p+1=12$ and $\langle \sigma, \tau \rangle\simeq M_{12}$, $M_{11}$ or
    $L_2(11)$;
    \item \label{FW:ix} $m=p+1=24$ and $\langle \sigma, \tau \rangle\simeq
    M_{24}$;
    \item \label{FW:x} $m=p+1$, $p$ is a Mersenne prime and $\langle \sigma, \tau \rangle$
    is a Frobenius group with kernel an elementary abelian
    $2$-group of order $m$ and complement of order $p$;
    \item \label{FW:xi} $m=p+1$, $p$ is a Mersenne prime, $H=L_k(2)$, $2^k=m$, $k\neq 3$ and $\langle \sigma, \tau \rangle$
    is isomorphic to the semi-direct product of an elementary
abelian $2$-group by $H$ with $H$ acting in its natural
action; 
    \item \label{FW:xii} $m=p+1=3$ and $\langle \sigma, \tau \rangle\simeq
    \s_3$; or
    \item \label{FW:xiii} $\langle \sigma, \tau \rangle\simeq \A_m$.
\end{enumerate}
\end{theorem}

The proof uses the list of 2, 3-transitive simple groups appearing in \cite{Ca}.

\subsection{On projective special linear groups}\label{subsec:PSL}
Let $r$, $k\in\N$, with $r$ a prime power. The projective special linear group $L_k(r)$ has order
\begin{align*}
|L_k(r)|= \frac{r^\frac{k(k-1)}{2}}{d} \prod_{i=2}^k (r^i-1),
\end{align*}
where $d=\gcd(k,r-1)$, see \cite{Ar}. Assume that $p:=\frac{r^k-1}{r-1}$ is prime. Then $k$ is prime and $d=1$. The group $L_k(r)$, which coincides with $\SL_k(r)$ in this case, is a primitive group contained in $\A_p$. Indeed, $L_k(r)$ is a 2-transitive permutation group of degree $p$, see \cite[Lemma 8.29]{Is}.
Notice that $\SL_k(r)$ has elements of order $p$ see \cite[Corollary 3]{Da}.
On the other hand, a Sylow $p$-subgroup of $L_k(r)$ has order $p$, it is self-centralizer and its normalizer has order $pk$.
Hence, the number of conjugacy classes of elements of order $p$ in $L_k(r)$ is even and equal to $(p-1)/k$. For more information on the number of conjugacy classes in finite classical groups see \cite{M} and \cite{W}.

\begin{obs}\label{rem:c:ccp>4}
Let $p$ be prime as above, with $p\geq 13$.

(a) We claim that $t:=(p-1)/k\geq 4$. Indeed, this is easy to see for $k=2$ and $k=3$. For $k\geq 5$, the result follows from $p\geq 2^k-1\geq 4k+1$ which can be proven by induction on $k$.

Let $\sigma$ be an element of order $p$ in $L_k(r)\subset \A_p$. We recall that $\oc\cap \la\sigma\ra=\{\sigma^\ell\,| J(\ell,p)=1\,\}$, where $J(\ell,p)$ is the Jacobi symbol of $\ell$ and $p$, see \cite[Claim 1, p. 240]{AFGV-alt}.
Let $\mathfrak K:=\{\C_1,\dots,\C_t\}$ be the set of conjugacy classes of elements of order $p$ in $L_k(r)$ and take $\ell$ such that $J(\ell,p)=-1$. The set $\mathfrak K$ splits into two sets $\mathfrak K_1:=\{\C_i\,|\,\C_i\subset \oc_\sigma^{\A_p}\}$ and $\mathfrak K_2:=\{\C_i\,|\,\C_i\subset \oc_{(1\, 2)\trid \sigma}^{\A_p}\}$.
It is easy to see that $\mathfrak K_1$ and $\mathfrak K_2$ have the same cardinality. Indeed, if $\C_i \in\mathfrak K_1$, then $\C_i^\ell:=\{x^\ell\,|\,x \in \C_i\} \in\mathfrak K_2$;
now, the result follows using that the function $\varphi_\ell:\oc_\sigma^{\A_p} \to \oc_{(1\, 2)\trid \sigma}^{\A_p}$, given by $\varphi_\ell(x)=x^\ell$, is bijective and it induces a bijection between $\mathfrak K_1$ and $\mathfrak K_2$.

\medbreak

(b) If $\oc_1$ and $\oc_2$ are two conjugacy classes of elements of order $p$ in $L_k(r)$, then there exist $\sigma\in \oc_1$ and $\tau\in \oc_2$ such that the order of the element $\sigma\tau$ is different from $1$, $2$ and $p$. This follows from \cite[Thm. 2]{Go} which states that
in any finite simple group of Lie type $G$ the product of any two conjugacy classes consisting of regular semisimple elements contains all semisimple elements of the group.
We recall that an element $g\in G$ is called \emph{regular semisimple} if the order of its centralizer in $G$ is relatively prime to the characteristic of the corresponding finite field. In our case, any element of order $p$ has centralizer of order $p$, thus it is regular semisimple since $p$ and $r$ are relatively prime.

\end{obs}

\subsection{}\label{sub:order<p^2}
Let $G$ be a finite group, $\oc$ a non-trivial conjugacy class of $G$, and $\sigma$, $\tau\in\oc$.
Assume that $(\sigma\tau)^2=(\tau\sigma)^2$; this amounts to saying that $\tau\sigma\tau$ commutes with $\sigma$ or, equivalently, that $\sigma\tau\sigma$ commutes with $\tau$.

\begin{lema}\label{le:}
If the centralizer of $\sigma$ in $G$ is cyclic of order $|\sigma|$, then the order of $\la\sigma,\tau\ra$ is at most $|\sigma|^2$.
If, in addition, $|\sigma|$ is prime, then $\sigma$ and $\tau$ commute or $|\sigma\tau|=2$.
\end{lema}
\pf
Since $\tau\sigma\tau$ commutes with $\sigma$, we have that $\tau\sigma\tau=\sigma^i$, for some $i$. Thus,
$\la\sigma,\tau\ra$ is at most $|\sigma|^2$. We also have that $\sigma\tau\sigma=\tau^j$, for some $j$. Then $\sigma^{i+1}=\tau^{j+1}$ 
because of the assumption $(\sigma\tau)^2=(\tau\sigma)^2$. Assume that $|\sigma|=p$, with $p$ prime. If $j+1 \neq 0 \mod (p)$, then $\tau\in \la\sigma\ra$, whereas if $j+1 = 0 \mod (p)$, then $\sigma\tau\sigma=\tau^{-1}$, and
$|\sigma\tau|=2$.
\epf

\begin{lema}\label{teor:usefulresult2}
Let $G$ be a finite group and let $\oc$ be a conjugacy class of $G$ whose elements have order $p$, with $p$ an odd prime. Assume that
\begin{enumerate}
\renewcommand{\theenumi}{\alph{enumi}}\renewcommand{\labelenumi}{(\theenumi)}
\item the centralizer in $G$ of an element in $\oc$ has order $p$, and
\item \label{eq:hipothesis3} there exists a subgroup $H$ of $G$ such that $\oc$ contains two different conjugacy classes $\oc_1$, $\oc_2$ of $H$.
\end{enumerate}
If for some $\sigma \in \oc_1$ fixed
\begin{align}\label{eq:inequality}
\text{there exists $\tau \in \oc_2$ such that  $|\sigma \tau|\neq 1, 2, p$},
\end{align}
then $\oc$ is of type D.
\end{lema}

\pf Let $\sigma\in\oc_1$. The condition \eqref{eq:inequality} implies
that there exists $\tau\in \oc_2$ such that $\tau$ does not commute with
$\sigma$ and $|\sigma \tau|\neq 2$. By the previous discussion,
$(\sigma\tau)^2\neq(\tau\sigma)^2$. Now, since $\sigma$ and $\tau$ are not
conjugated in $H$, the condition (Ax.~2) holds. \epf

\begin{cor}
Let $p$ be a prime number, $p\geq 5$, and let $\oc$ be the conjugacy class of $p$-cycles in $\s_p$. Then $\oc$ is of type D.
\end{cor}
\pf
It follows from Lemma \ref{teor:usefulresult2} with $\sigma=(1 \, 2\, \cdots \, p)$, $H=\A_p$, $\oc_1=\oc_\sigma^{\A_p}$, $\oc_2=\oc_{(1\, 2)\trid \sigma}^{\A_p}$ and
$\tau=(1\, 3)\trid \sigma$. 
\epf

\section{Proof of the main result}\label{se:proof}

Let $p\in \N$ be an odd prime, with $p\geq 5$, and $m\in\{p,p+1\}$.
Define the $p$-cycle $\sigma=(1 \, 2\, \cdots \, p)$ and let $\oc$ be the conjugacy class
of $\sigma$ in $\A_m$. We will determine when
there exists $\tau\in \oc$ such that
(Ax.~1) and (Ax.~2) hold using Theorem \ref{teor:FW}.
Notice that if $\la \sigma,\tau\ra\simeq \Z/p\Z$ or $\A_{m}$, with $\tau\in \oc$, then
(Ax.~1) does not hold for $\Z/p\Z$ and (Ax.~2) does not hold for $\A_m$.

\medbreak

(I) Assume that $m=p$.

Suppose that $p$ is not of the form $\frac{r^k-1}{r-1}$, with $r$ a prime power. If $p\neq 11, 23$, and $\tau \in \oc$, then $\la \sigma,\tau\ra\simeq \Z/p\Z$ or $\A_p$; hence, $\oc$ is not of type D. If $p=11$ and $\tau \in \oc$, then $H:=\la \sigma,\tau\ra\simeq \Z/11\Z$, $\A_{11}$, $L_2(11)$ or $M_{11}$.
In the last two cases, (Ax.~2) does not hold since each of the groups $L_2(11)$ and $M_{11}$ have two conjugacy classes of elements of order 11 and each of them is contained in different conjugacy classes in $\A_{11}$; indeed, if $h\in H$ and $|h|=11$, then each of the two conjugay classes of elements of order 11 in $\A_{11}$ contains
some power $h^\ell$, $1\leq \ell\leq 10$.
Hence, $\oc$ is not of type D. The case $p=23$ follows analogously.

Suppose that $p=\frac{r^k-1}{r-1}$, with $r$ a prime power. By Subsection \ref{subsec:PSL}, $\A_p$ contains a subgroup $H$ such that $\sigma\in H$ and $H\simeq L_k(r)$.
Assume that $p\geq 13$. By Remark \ref{rem:c:ccp>4} (a), there are at least two conjugacy classes $\oc_1$ and $\oc_2$ of $H$ contained in $\oc$. We can assume $\sigma\in \oc_1$; then condition (Ax.~2) holds for any $\tau\in \oc_2$.
By Lemma \ref{teor:usefulresult2} and Remark \ref{rem:c:ccp>4} (b), condition (Ax.~1) holds for some $\tau\in \oc_2$, and $\oc$ yields of type D.
Finally, if $p=5$ or $7$, then $\oc$ is not of type D. Indeed, if $H$ is a subgroup generated by two $p$-cycles, then $H\simeq \Z/5\Z$ or $\A_5\simeq L_2(4)$ when $p=5$, whereas $H\simeq \Z/7\Z$, $\A_7$ or $L_2(7)$ when $p=7$. Notice that $L_2(7)$ has only two conjugacy classes of elements of order $7$.

\medbreak

(II) Assume that $m=p+1$.

Suppose that $p$ is not of the form $\frac{r^k-1}{r-1}$, with $r$ a prime power. If $p\neq 11, 23$, and $\tau \in \oc$, then $\la \sigma,\tau\ra\simeq \Z/p\Z$, $\A_{p+1}$, $\A_{p}$ or $L_2(p)$. In the last two cases, condition (Ax.~2) does not hold since these groups have two conjugacy classes of elements of order $p$ (see \cite{FH} or \cite{Ad} for the groups $L_2(p)$) and each of them is contained in different conjugacy classes in $\A_{p+1}$. Hence, $\oc$ is not of type D.
If $p=11$ and $\tau \in \oc$, then $\la \sigma,\tau\ra\simeq \Z/11\Z$, $\A_{12}$, $\A_{11}$, $M_{11}$, $M_{12}$ or $L_2(11)$.
In the last four cases, condition (Ax.~2) does not hold since these groups have two conjugacy classes of elements of order 11 and each of them is contained in different conjugacy classes in $\A_{12}$. Hence, $\oc$ is not of type D.
If $p=23$ and $\tau \in \oc$, then $\la \sigma,\tau\ra\simeq \Z/23\Z$, $\A_{24}$, $\A_{23}$ or $M_{24}$; therefore, $\oc$ is not of type D as above.

Suppose that $p=\frac{r^k-1}{r-1}$, with $r$ a prime power. By (I) and Remark \ref{rem:OpOp+1} (b), if $p\geq 13$, then $\oc$ is of type D.
If $p=5$, then the subgroup generated by two $5$-cycles is $\Z/5\Z$, $\A_6$ or $\A_5$; thus, $\oc$ is not of type D.
Finally, if $p=7$, or more generally, $p=2^h-1$, with $h\geq 3$, is a Mersenne prime, then $\oc$ is of type D. Indeed, set $q=2^h$; then, the group $H = \F_{\hspace{-2pt}q}
\rtimes \F_{\hspace{-2pt}q}^{\times}\simeq \F_{\hspace{-2pt}q}
\rtimes \Z/p\Z$ acts on $\F_q$ by translations and dilations. Now, if we
identify $\{1, \dots, q\}$ with $\F_q$, then $H$ is isomorphic to
a subgroup of $\s_q$, see \cite{AFGV-alt}.

This finishes the proof of Theorem \ref{theorem:main}.

\begin{obs}\label{obs:ciclot:primos:1000}
The prime numbers $p$ of the form $(r^k-1)/(r-1)$, with $r$ a prime power and
$p<1000$ are:
$3=2^2-1$, 
$5=\frac{4^2-1}{4-1}$, 
$7=2^3-1$,
$13=\frac{3^3-1}{3-1}$,
$17=\frac{16^2-1}{16-1}$, 
$31=2^5-1=\frac{5^3-1}{5-1}$,
$73=\frac{8^3-1}{8-1}$,
$127=2^7-1$,
$257=\frac{256^2-1}{256-1}$, 
$307=\frac{17^3-1}{17-1}$,
$757=\frac{27^3-1}{27-1}$.

It is not known if the family of this kind of primes is finite or not. Indeed, it contains the families of Mersenne primes and Fermat primes. A discussion on numbers of this form can be found in \cite{EGSS}.

\end{obs}

\begin{obs}

(a)
The abelian subracks $T$ of $\oc$, with $\sigma\in T$, are contained in $\oc\cap \la\sigma\ra$, see Remark \ref{rem:c:ccp>4} (a). Thus, any maximal abelian subrack of $\oc$ has $(p-1)/2$ elements and it is isomorphic to $\oc\cap \la\sigma\ra$.

\medbreak

(b) Let $\oc_{(p)}$ be a conjugacy class of $p$-cycles in $\A_p$ not of type D. By the Theorem \ref{theorem:main}, $p=5$, $7$ or $p$ is not of the form $(r^k-1)/(r-1)$, with $r$ a prime power. It is clear that a subrack $X$ of $\oc_{(p)}$ is the union of the conjugacy classes $\oc_x^H$, $x\in X$, where $H$ is the subgroup of $\A_p$ generated by the elements of $X$.
Notice that $H$ is a simple group since it is generated by $p$-cicles in $\A_p$.

Clearly, $H$ is abelian if and only if $X$ is an abelian subrack. Assume that $H$ is not abelian. Then $H$ must be a $2$-transitive simple group of prime degree; this follows as in the step 5) of the proof of Theorem \ref{teor:FW} given in \cite{FW}. Then it occurs that $H$ is as in the cases (iii), (vi), (vii) or (xiii) of Theorem \ref{teor:FW}.
Hence, the non-abelian subracks of $\oc_{(p)}$ are conjugacy classes of elements of order $p$ in the subgroups appearing in that cases.

Therefore, the only cases where $\oc_{(p)}$ has proper non-abelian subracks are $p=7$, $11$ and $23$, and these subracks are
isomorphic to a conjugacy class of elements of order $p$ in $L_2(7)$, $L_2(11)$ or $M_{11}$, and $M_{23}$, respectively.
For instance, any proper non-abelian subrack $X$ of $\oc_{(p)}$ has $24$ elements for $p=7$ and $60$ or $720$ elements for $p=11$;
moreover, $X$ is not fixed by conjugation of any element in $\oc_{(p)}\setminus X$.

\medbreak

(c) Let $\oc$ be a conjugacy class of $p$-cycles in $\A_{p+1}$ not of type D. By the Theorem \ref{theorem:main}, $p=5$ or $p$ is not of the form $(r^k-1)/(r-1)$, with $r$ a prime power.
As in (b) above, the proper non-abelian subracks of $\oc$ are conjugacy classes of elements of order $p$ in the corresponding subgroups appearing in the proof of Theorem \ref{theorem:main}.
They are: $\A_5$ for $p=5$; $\A_{11}$, $L_2(11)$, $M_{11}$ and $M_{12}$ for $p=11$; $\A_{23}$, $M_{23}$ and $M_{24}$ for $p=23$; $\A_{p}$ and $L_2(p)$ otherwise.

\medbreak

(d)  For the racks $\oc$ described in (b) and (c) above it would be possible to decide if the dimension of $\toba(\oc,q)$ is not finite, for some $2$-cocycle $q$, as mentioned in
\cite[\S 2.6]{AFGaV-cla}. Indeed, for an abelian subrack $T$ of $\oc$ and any $2$-cocycle $q$ we could determine if the diagonal braiding associated with $q$ gives rise to a Nichols algebra of infinite dimension; in that case $\toba(\oc,q)$ would be also of infinite dimension. For this we can use the classification of finite-dimensional Nichols algebras of diagonal type \cite{H}.

In that sense, we compute using \textsf{RiG} (see \cite{GV}) the second abelian rack cohomology group of some of this racks:
\begin{itemize}
  \item $H^2(\oc_{(5)},\ku^\times)=\ku^\times \times \G_{10}$. Notice that $\oc_{(5)}\simeq Q_{12,3}$ in \cite{V}.
  \item $H^2(X,\ku^\times)=\ku^\times \times \G_{14}$, with $X\subset\oc_{(7)}$, $|X|=24$.
\end{itemize}
See \cite{AG} for the considered cohomology theory of racks and \cite{AFGaV-cla} for the use of \textsf{RiG} in our cases.
It would be expected that $H^2(\oc_{(p)},\ku^\times)=\ku^\times \times \G_{2p}$.
\end{obs}

\subsection*{Acknowledgement}
The author is grateful to the referee for several helpful suggestions.
I thank R. Guralnick and A. Hulpke for very important comments on linear simple groups.
I also thank L. Vendramin for interesting discussions about \textsf{RiG}.
Part of this work was done during a postdoctoral position at Universit\'e Paris
Diderot -- Paris 7; I am very grateful to Prof. Marc Rosso for his kind hospitality.

\end{document}